\documentclass[final]{siamltex}
\usepackage{amssymb, amsmath}
\usepackage{float,epsfig}
\usepackage{algpseudocode}
\usepackage{algorithm}

\newtheorem{remark}[theorem]{ Remark}
\newtheorem{exam}[theorem]{\bf Example}

\newcommand{\ba}{\begin{array}}
\newcommand{\ea}{\end{array}}

\newcommand{\be}{\begin{equation}}
\newcommand{\ee}{\end{equation}}
\newcommand{\beano}{\begin{eqnarray*}}
\newcommand{\eeano}{\end{eqnarray*}}

\def\R{{\mathbb R}}
\def\C{{\mathbb C}}

\def\lam{\lambda}
\def\sig{\sigma}

\def\diag{\mathrm{diag}}

\def\rank{\mathrm{rank}}

\def \nrank{\mathrm{nrk}}

\def \sp{\mathrm{Sp}}

\title{Solution Method for Higher Order System }

\author{ Namita Behera \thanks{Department of Mathematics, Sikkim University,
Sikkim-737102, INDIA,({\tt nbehera@cus.ac.in, niku.namita@gmail.com })}}

\begin{document}

\maketitle

\begin{center}
	\today
\end{center}

\begin{abstract}
Consider a higher order state space system and  the aim of this paper is to linearize the system preserving system characteristics. That is, linearization preserving system characteristics(e.g, controllability, observability, various zeros  and transfer function) for analysis of higher order systems gives the solution for higher order system.  We study recovery of zero directions of higher order state space system from those of the linearizations. That is, the zero directions of the transfer functions associated to higher order state space system are recovered from the eigenvectors of the Fiedler pencils without performing any arithmetic operations.
\end{abstract}

\begin{keywords}
Higher order system, System matrix, zero direction, matrix polynomial, eigenvalue,  eigenvector,  matrix pencil, linearization, Fiedler pencil.
\end{keywords}

\begin{AMS}
65F15, 15A18, 15B57, 15A22
\end{AMS}

\section{Introduction}
 We denote by $\C[\lam]$ the polynomial ring over the complex field $\C.$  Further, we denote by $ \C^{m\times n}$ and $  \C[\lam]^{m\times n}$, respectively, the vector spaces of $m\times n$ matrices and  matrix polynomials over $\C.$  A matrix polynomial $P(\lam)$ is said to be regular if $\det(P(\lam)) \neq 0$ for some $\lam \in \C$.
Throughout this paper, we consider the higher order linear time invariant (LTI) state space system $\Sigma_1$ given by
\begin{equation*}
\begin{aligned}
P\left(\frac{d}{dt}\right) x(t) &= B u(t), \\
y(t) &= C x(t) + D u(t),
\end{aligned}
\end{equation*}
where $P(\lam) = \sum_{j=0}^{m}\lam^{j}A_j \in \C[\lam]^{n \times n}$ be regular and $ D \in \C^{r \times r}, C \in \C^{r \times n}, B \in \C^{n \times r}$.

Higher order state space system plays an important role in system theory. We wish to compute zeros of higher order systems. So, consider a higher order state space system $\Sigma_1$ and derive linearized state space systems which are strict system equivalent to the higher order system.

\label{sst}
Consider a linear, time invariant (LTI), multivariable system $\Sigma$ given by
\begin{eqnarray*}
 A\left(\frac{d}{dt}\right) x(t) &=& B\left(\frac{d}{dt}\right) u(t) \nonumber \\
 y(t) &=& C\left(\frac{d}{dt}\right) x(t) + D\left(\frac{d}{dt}\right) u(t), \,\, t \geq 0
\end{eqnarray*}
where $\frac{d}{dt}$ is the differential operator, $ A(\lam) \in \mathbb{C}[\lam]^{r \times r}$ is regular, $B(\lam) \in \mathbb{C}[\lam]^{r \times m},$ $ C(\lam) \in \mathbb{C}[\lam]^{p \times r}, D(\lam) \in \mathbb{C}[\lam]^{p \times m}, u(t) : \R^{+} \rightarrow \R^{m}$ is the input vector, $x(t) : \R^{+} \rightarrow \R^{r}$ is the state vector and $y(t) : \R^{+} \rightarrow \R^{p}$ is the output vector \cite{vardulakis}.

\begin{definition}\cite{rosenbrock70}
The matrix polynomial $\mathcal{S}(\lam)$ given by $$\mathcal{S}(\lam ) = \left[
                            \begin{array}{c|c}
                              -A(\lam) & B(\lam) \\
                              \hline
                              C(\lam) & D(\lam) \\
                            \end{array}
                          \right] \in \C[\lam]^{(r+p)\times (r+m)} $$ is called the Rosenbrock system matrix or the  Rosenbrock system polynomial or simply called system matrix of the system $\Sigma$. The rational matrix $G(\lam)$ given by $$G(\lam) = D(\lam) + C(\lam) A(\lam)^{-1} B(\lam)  \in \C(\lam)^{p \times m}$$ is called the transfer function of the system $\Sigma$.
\end{definition}


Now, we show that the state-space framework so developed in \cite{rafinami} could be gainfully used to linearize a higher order LTI state-space system so as to analyze and solve the higher order LTI system. That is, linearization preserving system characteristics(e.g, controllability, observability, various zeros  and transfer function) for analysis of higher order systems gives the solution for higher order system. Furthermore, we show that the linearized systems are strict system equivalent to the higher order systems and hence preserve system characteristics of the original systems.

The paper is organized as follows. Section $2$ presents some basic concepts on system theory and index tuples which we need throughout this paper.  In section $3$ we show that the state-space framework so developed in \cite{rafinami} could be gainfully used to linearize (Fiedler-like linearizations) a higher order LTI state-space system so as to analyze and solve the higher order LTI system. Further, in the same section we show that the linearized systems are strict system equivalent to the higher order systems and hence preserve system characteristics of the original systems. Furthermore, we study recovery of zero directions of $\Sigma_1$ from those of the Fiedler-like pencils of $G(\lam)$.

{\bf Notation.} An $m\times n$ rational matrix function $G(\lam)$ is an $m\times n$ matrix whose entries are rational functions of the form $p(\lam)/{q(\lam)},$ where $p(\lam)$ and $q(\lam)$ are scalar polynomials in $\C[\lam].$  An $n\times n$ rational matrix function $G(\lam)$ is said to be {\em regular} if  $\rank(G(\lam)) = n$ for some  $\lam \in \C.$ If $G(\lam)$ is regular then $\mu \in \C$ is said to be an {\em eigenvalue} of $G(\lam)$ if $ \rank(G(\mu)) < n.$ We denote  the $j$-th column of the $n \times n$ identity matrix $I_n$ by $e_j$ and the transpose of a matrix $A$ by $A^T.$ We denote the Kronecker product of matrices $A$ and $B$ by $A \otimes B.$

\section{Basic Concepts}
The normal rank~\cite{rosenbrock70} of $ X(\lam) \in \C(\lam)^{m\times n}$ is denoted by $ \nrank(X)$ and is given by $\nrank(X) := \max_{\lam}\rank(X(\lam))$ where the maximum is taken over all $\lam \in \C$ which are not poles of the entries of $X(\lam).$ If $\nrank(X) =n = m$ then $X(\lam)$ is said to be regular, otherwise $X(\lam)$ is said to be singular.

A complex number $\lam$ is said to be an eigenvalue of the system matrix $\mathcal{S}(\lam)$ if
$\rank(\mathcal{S}(\lam)) < \nrank(\mathcal{S}).$ An eigenvalue $\lam$ of $\mathcal{S}(\lam)$ is called an invariant zero
of the system $\Sigma$. We denote the set of eigenvalues of $\mathcal{S}(\lam)$ by $\sp(\mathcal{S}).$

Let $G(\lam)\in \mathbb{C}(\lam)^{m \times n}$ be a rational matrix function and let
\be \label{smform}
\mathbf{SM}(G(\lam)) = \diag\left( \frac{\phi_{1}(\lam)}{\psi_{1}(\lam)}, \cdots, \frac{\phi_{k}(\lam)}{\psi_{k}(\lam)}, 0_{m-k, n-k}\right)\ee
be the Smith-McMillan form~\cite{kailath, rosenbrock70} of  $G(\lam),$
where the scalar polynomials $\phi_{i}(\lam)$ and $ \psi_{i}(\lam)$ are monic, are pairwise coprime and,  $\phi_{i}(\lam)$ divides $\phi_{i+1}(\lam)$ and $\psi_{i+1}(\lam)$ divides $\psi_{i}(\lam),$ for $i= 1, 2, \ldots, k-1$.
The polynomials $\phi_1(\lam), \ldots, \phi_k(\lam)$ and $ \psi_1(\lam), \ldots, \psi_k(\lam)$ are called {\em invariant zero polynomials} and {\em invariant pole polynomials} of $G(\lam),$ respectively. Define
$\phi_{G}(\lam) := \prod _{j=1}^{k} \phi_{j}(\lam) \,\,\, \mbox{ and } \,\,\, \psi_{G}(\lam) := \prod _{j=1}^{k} \psi_{j}(\lam).$
 A complex number $ \lam $ is said to be a  zero of $G(\lam)$ if $ \phi_G(\lam) =0$ and
a complex number $ \lam$ is said to be a  pole of $G(\lam)$ if $\psi_G(\lam) =0.$ The {\bf spectrum} $\sp(G)$ of $G(\lam)$ is given by  $ \sp(G) :=\{ \lam \in \C : \phi_G(\lam) = 0\}.$ That is $\sp(G)$ is the set of zeros of $G(\lam),$ see \cite{rafinami}.

The right and the left null spaces of $G(\lam)$ are given by $ \mathcal{N}_r(G) := \{ x(\lam) \in \C(\lam)^n: G(\lam)x(\lam) = 0\}$ and $\mathcal{N}_l(G) :=\{ y(\lam) \in \C(\lam)^m : y^T(\lam) G(\lam) =0 \}.$

\begin{definition}\cite{mmmm06}
Matrix polynomial $U(\lam)$ is said to be unimodular if $\det U(\lam)$ is a nonzero constant, independent of $\lam$. Two matrix
polynomials $P(\lam)$ and $Q(\lam)$ are said to be equivalent if there exists unimodular matrices $U(\lam)$ and $V(\lam)$, such that $Q(\lam) = U(\lam)P(\lam)V (\lam).$ If $U(\lam), V (\lam)$ are constant matrices, then $P(\lam)$ and $Q(\lam)$ are said to be strictly equivalent.
\end{definition}

Let $P(\lam), Q(\lam) \in \C[\lam]^{m \times n}.$ Then $Z(\lam) \in \C^{m \times m}$ is said to be a {\bf common left divisor} of $P(\lam)$ and $Q(\lam)$ if $P(\lam) = Z(\lam)R(\lam)$ and $Q(\lam) = Z(\lam) T(\lam)$, for some $R(\lam), T(\lam) \in \C[\lam]^{m\times n}$. Equivalently, $\left[
\begin{array}{cc}
P(\lam) & Q(\lam) \\
\end{array}
\right] = Z(\lam)\left[
                   \begin{array}{cc}
                     R(\lam) & T(\lam) \\
                   \end{array}
                 \right].
$ Further, $Z(\lam)$ is said to be the {\bf greatest common left divisor} of $P(\lam)$ and $Q(\lam)$ if $U(\lam)$ is a common left divisor of $P(\lam)$ and $Q(\lam)$ then $Z(\lam) = U(\lam) V(\lam)$ for some $V(\lam) \in \C^{m \times m}$, see \cite{vardulakis}. \\

The matrix polynomial $P(\lam)$ and $Q(\lam)$ are said to be {\bf left coprime} if the greatest common left divisors of $P(\lam)$ and $Q(\lam)$ is unimodular. The greatest common right divisor of $P(\lam)$ and $Q(\lam)$  is defined similarly. Then $P(\lam)$ and $Q(\lam)$ are said to be {\bf right coprime} if the greatest common right divisor of $P(\lam)$ and $Q(\lam)$ is unimodular, see \cite{vardulakis}.

\begin{definition} \cite{vardulakis}
The LTI system $\Sigma$ is said to be of least order (or minimal) if $A(\lam)$ and $B(\lam)$ are left coprime, and $A(\lam)$ and $C(\lam)$ are right coprime.
\end{definition}

If $A(\lam) = \lam I -A, B(\lam) = B$ and $C(\lam) = C$ then:
\begin{itemize}
\item [(i)] $A(\lam)$ and $B(\lam)$ are left coprime $\Leftrightarrow$ $\rank [\lam I-A, B] = n$ for $\lam \in \C$,

\item [(ii)] $A(\lam)$ and $C(\lam)$ are right coprime $\Leftrightarrow$ $\rank \left[
                                                                                  \begin{array}{c}
                                                                                    \lam I-A \\
                                                                                    C \\
                                                                                  \end{array}
                                                                                \right] = n
$ for $\lam \in \C$.
 \end{itemize}

\begin{definition}\cite{vardulakis}
The degree of a matrix polynomial  $P(\lam) \in \C[\lam]^{p\times m}$, denoted by $\deg P(\lam)$, is defined as the maximum degree among the degrees of all its maximum order (non-zero) minors.
\end{definition}

Note that if $P(\lam)$ is square and regular, then $\deg P(\lam) = \det(P(\lam)).$

\begin{definition}\cite{vardulakis}
Let $\mathcal{S}_1(\lam)$ and $\mathcal{S}_2(\lam)$ be system matrices with extended forms
\begin{equation}
\mathcal{S}_{i}^{e}(\lam) = \left[
                \begin{array}{cc|c}
                    I_{q-r_i} & 0 & 0 \\
                    0 & A_{i}(\lam) & B_{i}(\lam) \\
                    \hline
                    0 & C_{i}(\lam) & D_{i}(\lam) \\
                   \end{array}
                   \right] \indent i = 1, 2,
\end{equation}
where $ A_i(\lam) \in \mathbb{C}[\lam]^{r_i \times r_i}$, $B_i(\lam) \in \mathbb{C}[\lam]^{r_i \times m}, C(\lam) \in \mathbb{C}[\lam]^{p \times r_i}, D(\lam) \in \mathbb{C}[\lam]^{p \times m},$ $n_i = \deg |(A_i(\lam))|$ and $ q \geq \max\{n_1, n_2\}$. Then $\mathcal{S}_1(\lam)$ and $\mathcal{S}_2(\lam)$ are said to be {\em strict system equivalent} (SSE) if there exists unimodular $U(\lam), V(\lam) \in \C[\lam]^{q \times q}$, and $X(\lam) \in \C[\lam]^{p \times q}$,  $Y(\lam) \in \C[\lam]^{q \times m}$ such that
$$\left[
    \begin{array}{cc}
      U(\lam) & 0 \\
      X(\lam) & I_p \\
    \end{array}
  \right] \left[
                \begin{array}{ccc}
                    I_{q-r_1} & 0 & 0 \\
                    0 & A_{1}(\lam) & B_{1}(\lam) \\
                    0 & C_{1}(\lam) & D_{1}(\lam) \\
                   \end{array}
                   \right] = \left[
                \begin{array}{ccc}
                    I_{q-r_2} & 0 & 0 \\
                    0 & A_{2}(\lam) & B_{2}(\lam) \\
                    0 & C_{2}(\lam) & D_{2}(\lam) \\
                   \end{array}
                   \right] \left[
                             \begin{array}{cc}
                               V(\lam) & Y(\lam) \\
                               0 & I_m \\
                             \end{array}
                           \right].
$$
\end{definition}

\begin{theorem}\cite{vardulakis} \label{vsset}
Every extended system matrix $\mathcal{S}_{i}^{e}(\lam) \in \C[\lam]^{(n+p)\times (n+m)}$ is a strict system equivalent to a system matrix in state space form.
\end{theorem}

\begin{definition}\cite{MFT}
An ordered tuple of indices consisting of consecutive integers is called a {\em string} and denoted by $(t:p)$ for the string of integers from $t$ to $p$, i.e.,
$$(t:p) := \begin{cases}
(t, t+1, \ldots, p), & \text{  if } t \leq p \\
\phi,  & \text{  if } t > p.
\end{cases}$$
\end{definition}

\begin{remark}
In the above definition, if $t_1 > p$ and $t_2 > p$, then both index tuples $(t_1 : p)$ and $(t_2 :p)$ correspond to the empty index tuple. To avoid this notation, we will adopt the notation $(\infty: p)$ for any tuple of the form $(t:p)$ having $t>p$ where applicable.
\end{remark}

\begin{definition} \cite{VAEN}
Let $h$ be a non-negative integer and $\textbf{q}$ be an index tuple containing indices from $\{0, 1, \ldots, h\}$. Then $\textbf{q}$ is said to be in {\em column standard form} if $$\textbf{q} = (a_{p}: b_{p}, a_{p-1}: b_{p-1}, \ldots, a_{2}: b_{2}, a_{1}: b_{1}), $$ with $0 \leq b_{1} < b_{2} < \ldots < b_{p-1} <b_{p} \leq h$ and $0 \leq a_{j} \leq b_{j}$, for all $j = 1, \ldots, p$. We denote this tuple by $csf(\textbf{q})$.
\end{definition}

\begin{lemma}\cite{MFT}
Let $\sig$ be a permutation of $\{h_0, h_{0} + 1, \ldots, h \},$ with $0 \leq h_0 \leq h.$ Then $\sig$ is in column standard form if and only if
$$\sig = (t_{\alpha-1} + 1 : h, t_{\alpha-2} + 1 : t_{\alpha-1}, \ldots, t_2 + 1 : t_3, t_1 + 1 : t_2, h_0 : t_1), $$
for some positive integers $h_0 \leq t_1 < t_2 < \cdots < t_{\alpha-1} < h.$

Denote $t_0 = h_0-1$ and $t_{\alpha} = h.$ We call each sequence of consecutive integers $(t_{i-1}+1 :t_i), $ for $i = 1, \ldots, \alpha,$ a string in $\sig$.
\end{lemma}

\begin{definition}\label{dci} \cite{TDM, rafinami3}
Let $\sigma : \{0, 1, \ldots, m-1\} \rightarrow \{1, 2, \ldots, m\}$ be a bijection.
\begin{itemize}

\item[(1)] For $d = 0, \ldots, m-2$, we say that $\sigma$ has a consecution at $d$ if $\sigma(d) < \sigma(d+1)$ and $\sigma$ has an inversion at $d$ if $\sigma(d) > \sigma(d+1)$.

\item[(2)] We denote the total number of consecutions in $\sigma$ by $c(\sigma)$ and the total number of inversions in $\sigma$ by $i(\sigma)$.

\item[(3)] The consecution-inversion structure sequence of $\sigma$ is the tuple $(c_{1}, i_{1}, c_{2}, i_{2}, \ldots, c_{l}, i_{l})$, where $\sigma$ has $c_{1}$ consecutive consecutions at $0, 1, \ldots, c_{1}-1;$ $i_{1}$ consecutive inversions at $c_{1}, c_{1}+1, \ldots, c_{1}+i_{1}-1$ and so on, up to $i_{l}$ inversions at $m-1-i_{l}, \ldots, m-2$, and is denoted by CISS$(\sigma)$.

\item[(4)] If $0 \in \sig$, then we refer to CIP$(\sig) := (c_0, i_0)$ as the consecution-inversion pair
of $\sig$ at $0$, where $c_0$ (resp., $i_0$) is the number of consecutions (resp., inversions)
of $\sig$ at $0$.
\end{itemize}
\end{definition}

Observe that if CIP$(\sig)= (c_0, i_0)$, then either CIP$(\sig) = (0, 0)$ or CIP$(\sig) = (c_0, 0)$
with $c_0 > 0$ or CIP$(\sig) = (0, i_0)$ with $i_0 > 0$.

\begin{definition}~\cite{TDM}
Let $P(\lambda) = A_{0} + \lambda A_{1} + \cdots + \lambda^{m}A_{m}$ be a matrix polynomial of degree $m$. For $k = 0, \ldots m$, the degree $k$ Horner shift of $P(\lambda)$ is the matrix polynomial $P_{k}(\lambda) := A_{m-k} + \lambda A_{m-k+1} + \cdots + \lambda^{k}A_{m}$. These Horner shifts satisfy the following:
\begin{align}
& P_{0}(\lambda) = A_{m}, \nonumber  \,\,\,
 P_{k+1}(\lambda) = \lambda P_{k}(\lambda) + A_{m-k-1} , \mbox{  for } 0\leq k\leq m-1, \,\,\,
P_{m}(\lambda) = P(\lambda). \nonumber \label{hs}
\end{align}
\end{definition}

\begin{definition}
Let $\mathcal{H} = (\mathcal{H}_{ij})$ be a block $m \times n$ matrix with $p \times q$ blocks $\mathcal{H}_{ij}$. The {\em block transpose} of $\mathcal{H}$ is the block $n \times m$ matrix with $p \times q$ blocks, and defined by $\mathcal{H}^{\mathcal{B}} = (\mathcal{H}_{ji})_{n \times m}.$
\end{definition}

\begin{theorem}\cite{behera}\label{ralitfmp}
Suppose that $P(\lam)$ is regular. Let $L_{\sig}(\lam)$ be the linearization of $P(\lam)$ associated with a bijection $\sig : \{0, 1, \ldots, m-1\} \rightarrow \{1, 2, \ldots, m\}. $ Let $U(\lam)$ and $V(\lam)$ be unimodular matrices such that $$U(\lam)L_{\sig}(\lam) V(\lam) = \left[
                                                                                                                        \begin{array}{c|c}
                                                                                                                          I_{(m-1)n} &  \\
                                                                                                                          \hline
                                                                                                                           & P(\lam) \\
                                                                                                                        \end{array}
                                                                                                                      \right].
 $$
\begin{itemize}

\item [(a)] Define $E_{\sig}(P) : \C^{n} \rightarrow \C^{mn}$ by $E_{\sig}(P) = V(\lam) (e_m \otimes I_n) $ and $F_{\sig}(P) : \C^{mn} \rightarrow \C^{n} $ by $F_{\sig}(P) = (e_m^{T}\otimes I_n)V(\lam)^{-1}$. Then $F_{\sig}(P)E_{\sig}(P) = I_n$.

\item[(b)] Further, if CISS$(\sig) = (c_1, i_1, \ldots, c_l, i_l)$. Then
\begin{equation*}
F_{\sig}(P) = e_{(m-c_1)}^{T} \otimes I_n.
\end{equation*}

\item [(c)] \cite{TDM} Furthermore, let $P_0, \ldots, P_m$ be the  Horner shifts of $P(\lambda)$ and $csf(\sigma) = (\textbf{b}_{\beta}, \ldots, \textbf{b}_{1})$, where $\textbf{b}_{k} = (t_{k-1}+1:t_{k})$, for $k = 1, \ldots, \beta $. Then
\begin{align}
E_{\sigma}(P) = \left[
      \begin{array}{c}
        B_{0}\, B_{ 1}\, \ldots\, B_{m-1} \\
      \end{array}
    \right]^{\mathcal{B}}   \label{rqoevfp}
\end{align}
where, if $\sigma(i) \in \textbf{b}_{k}$, for some $k = 1, \ldots, \beta $, then
\begin{equation*}
B_{i} = \begin{cases}
\lambda^{k-1}I_n, & \text{if } i = m-t_{k}-1 \\
\lambda^{k-1}P_{i}, & \text{otherwise. }
\end{cases} \label{biforev}
\end{equation*}

\end{itemize}
Thus $E_{\sig}(P) : \mathcal{N}_r(P(\lam)) \rightarrow \mathcal{N}_r(L_{\sig}(\lam))$ and $F_{\sig}(P) : \mathcal{N}_r(L_{\sig}(\lam)) \rightarrow \mathcal{N}_r(P(\lam))$ are isomorphisms.
\end{theorem}

The reversal of a bijection plays an important role in showing that there is an isomorphism between the left nullspace of $\mathcal{S}(\lam)$ and the left nullspace of pencil $\mathbb{L}_{\sigma}(\lam)$.

\begin{definition} \cite{TDM}
The reversal rev$\,\sigma$ of a bijection $\sigma : \{0, 1, \ldots, m-1\}\rightarrow \{1, 2, \ldots, m\}$ is a bijection from $\{0, 1, \ldots, m-1\}$ into $\{1, 2, \ldots, m\}$ defined by rev$\,\sigma(i) = m+1-\sigma(i)$ for $0 \leq i \leq m-1$. Equivalently, if $\sig = (\sig(0), \ldots, \sig(m-1))$ then rev $\sig = (\sig(m-1), \ldots, \sig(0)). $
\end{definition}

%
%
%

For the left nullspace of $P(\lam)$ and $L_{\sig}(\lam)$, we have the following.

\begin{theorem}\cite{rafinami2}\label{dsigp}
Let $L_{\sig}(P):= \lam M_m - M_{\sig}$ be the Fiedler linearization of $P(\lam)$ associated with a bijection $\sig$. Let $E_{\sig}(P)$ and $F_{\sig}(P)$ be as in Theorem \ref{ralitfmp}. Define $H_{\sig}(P): \C^{n} \rightarrow \C^{mn}$ and $K_{\sig}(P) : \C^{mn} \rightarrow \C^{n} $ by $H_{\sig}(P) = U(\lam)^{T}(e_m \otimes I_n)$ and $K_{\sig}(P) = (e_m^{T} \otimes I_n)(U(\lam)^{-1})^{T}$. Then $K_{\sig}(P)H_{\sig}(P) = I_n. $ Further, $H_{\sig}(P) = E_{rev \, \sig}(P^{T})$ and  $K_{\sig}(P) = F_{rev \,\sig}(P^{T}). $ Furthermore, if CISS$(\sig) = (c_1, i_1, \ldots, c_l, i_l)$, then
$$K_{\sig}(P) =
\begin{cases}
e_{m}^{T} \otimes I_n, \,\,\, \text{ if } c_1> 0 \\
e_{(m-i_1)}^{T} \otimes I_n, \,\,\, \text{ if } c_1=0.
\end{cases}$$
Thus $H_{\sig}(P): \mathcal{N}_{l}(P(\lam)) \rightarrow \mathcal{N}_{l}(L_{\sig}(\lam))$ and $K_{\sig}(P) : \mathcal{N}_l(L_{\sig}(\lam)) \rightarrow \mathcal{N}_l(P(\lam))$ are isomorphisms.

\end{theorem}

Now, we consider a higher order state space system and derive linearized state space systems which are strict system equivalent to the higher order system.

\section{Higher order state space system}
Now, we consider a higher order LTI system. We wish to compute zeros of higher order systems.
Let $P(\lam) = \sum\limits_{j=0}^{m}\lam^{j}A_j \in \C^{n \times n}[\lam] $ be regular. Consider the {\em higher order LTI state space system} $\Sigma_1$ given by
\begin{equation}
\begin{aligned}
P\left(\frac{d}{dt}\right) x(t) &= B u(t), \\
y(t) &= C x(t) + D u(t),
\end{aligned}\label{mulsss}
\end{equation}
where $D \in \C^{r \times r}, C \in \C^{r \times n}, B \in \C^{n \times r}$. The associated system matrix is given by
\begin{equation}
\mathcal{S}(\lam) = \left[
                       \begin{array}{c|c}
                         -P(\lam) & B \\
                         \hline
                         C & D \\
                       \end{array}
                     \right] \in \C^{(n+r)\times (n+r)} \label{smh}
\end{equation}
and the transfer function of $\Sigma_1$ is given by
\begin{equation}
G(\lam) = C P(\lam)^{-1} B + D \in \C^{r \times r}. \label{tfh}
\end{equation}

Consider matrix polynomial $P(\lam) = \sum_{i=0}^{m}\lam^{i} A_i, A_i \in \C^{n \times n}$,then Fiedler matrices associated with $P(\lam)$  are  $nm \times  nm$  matrices $ M_0,\ldots, M_m$  given by~\cite{TDM}
$$
 M_{0} := \left[
                   \begin{array}{cc}
                     I_{(m-1)n} &  \\
                      & -A_{0} \\
                   \end{array}
                 \right], \, M_{i} := \left[
  \begin{array}{cccc}
    I_{(m-i-1)n} &  & & \\
     & -A_{i} & I_{n} & \\
     & I_{n} & 0  &  \\
     &   &   & I_{(i-1)n}\\
  \end{array}
\right] $$
and
\begin{equation}
M_{m} := \left[
          \begin{array}{cc}
            A_{m} &  \\
             & I_{(m-1)n} \\
          \end{array}
        \right],
\,  M_{i}^{-1} = \left[
  \begin{array}{cccc}
    I_{(m-i-1)n} &  & & \\
     & 0 & I_{n} & \\
     & I_{n} & A_{i}  &  \\
     &   &   & I_{(i-1)n} \\
  \end{array}
\right], \label{ifp}
\end{equation}   $  i= 1, \ldots, m-1. $
The Fiedler matrices $M_i$ and $M_j$ commute when $|i-j| >1,$ that is,  $M_{i}M_{j} = M_{j}M_{i}$ when $|i - j| > 1.$ Let  $\sigma : \{0, 1, \ldots, m-1\}\rightarrow \{ 1, 2, \ldots, m\}$ be a bijection. Then $M_{\sigma}$ denotes the product of Fiedler matrices given by $M_{\sigma} := M_{\sigma^{-1}(1)}M_{\sigma^{-1}(2)}\cdots M_{\sigma^{-1}(m)}. $  Note that $\sigma(i)$ describes the position of the factor $M_{i}$ in the product $M_{\sigma}$, that is, $\sigma(i) = j$ means that $M_{i}$ is the $j$-th factor in the product. For convenience, we set   $M_{\emptyset} := I_{nm},$ where $\emptyset$ is the empty set. The $nm\times nm$ pencil $L_\sigma(\lam)$ given by  $L_{\sig}(\lam) := \lam M_m - M_\sig$ is called the Fiedler pencil of $P(\lam)$ associated with $\sig,$ see~\cite{TDM, AV}.

Now, we define the Fiedler pencil of $\mathcal{S}(\lam)$ given in (\ref{smh}). Define $(nm +r) \times (nm+ r)$ matrices
\begin{equation}
\mathbb{M}_{0} := \left[
                       \begin{array}{c|c}
                         M_{0} & -e_{m} \otimes B \\
                         \hline
                         -e_{m}^{T}\otimes C & -D \\
                       \end{array}
                     \right], \indent \mathbb{M}_{m} := \left[
                                        \begin{array}{c|c}
                                          M_{m} & 0 \\
                                          \hline
                                          0 & 0_r \\
                                        \end{array}
                                      \right] \label{0mmfrh}
\end{equation}
and
\begin{equation}
\mathbb{M}_{i} := \left[
                       \begin{array}{c|c}
                         M_{i} & 0 \\
                         \hline
                         0 & I_{r} \\
                       \end{array}
                     \right] \indent i = 1, \ldots, m-1, \label{imfrh}
\end{equation}
where $M_i, i=0:m$ are the Fiedler matrices associated with $P(\lambda)$. We refer $\mathbb{M}_i, i=0:m$ are the Fiedler matrices associated with $\mathcal{S}(\lam)$. The Fiedler matrices $\mathbb{M}_{i}$ associated with $\mathcal{S}(\lam)$ have the following properties.
\begin{itemize}

\item The matrix $\mathbb{M}_{m}$ is always singular.

\item Since $M_{i} M_{j} = M_{j}M_{i}$ for $|i - j| > 1$ so it follows that
\begin{equation}
\mathbb{M}_{i} \mathbb{M}_{j} = \mathbb{M}_{j} \mathbb{M}_{i} \indent\text{  for   } |i-j| > 1, \text{ except for } \mathbb{M}_{m} \text{   and  } \mathbb{M}_{0}.
\end{equation}

\item Note that since $M_i, i=1:m-1$ are invertible so, the Fiedler matrices $\mathbb{M}_{i}, i=1:m-1$ are invertible.
\end{itemize}

Let $\sigma : \{0, 1, \ldots, m-1\}\rightarrow \{ 1, 2, \ldots, m\}$ be a bijection. Then  $$\mathbb{M}_{\sigma}:=\mathbb{M}_{\sigma^{-1}(1)}\mathbb{M}_{\sigma^{-1}(2)}\cdots \mathbb{M}_{\sigma^{-1}(m)}.$$ Observe that if $\mathbb{M}_{\sigma} = \mathbb{M}_{\sigma^{-1}(1)}\mathbb{M}_{\sigma^{-1}(2)}\cdots \mathbb{M}_{\sigma^{-1}(m)}$ then
$$\mathbb{M}_{\textrm{rev}\,\sigma} = \mathbb{M}_{\sigma^{-1}(m)}\mathbb{M}_{\sigma^{-1}(m-1)}\cdots \mathbb{M}_{\sigma^{-1}(1)}. $$
If $\sigma$ and $\tau$ are two bijections from $\{0, 1, \ldots, m-1\}$  to $ \{ 1, 2, \ldots, m\}$ then because of commutation relation we may have $\mathbb{M}_\sigma = \mathbb{M}_\tau.$ Consequently, a Fiedler pencil may be associated with more than one  bijection.

Now, we define Fiedler pencil of $G(\lam)$ in (\ref{tfh}).

\begin{definition}[Fiedler pencil]\label{dfphos}
Let $\mathcal{S}(\lam)$ be the system matrix given in (\ref{smh}). For any bijection $\sigma : \{0, 1, \ldots, m-1\} \rightarrow \{1, 2, \ldots, m\}$, we define the pencil $\mathbb{S}_{\sig}(\lam)$ given by
$\mathbb{S}_{\sig}(\lam) = \lambda \mathbb{M}_{m}- \mathbb{M}_{\sigma} $ is said to be the Fiedler pencil of system polynomial $\mathcal{S}(\lam)$ associated with $\sig$. We also refer to $\mathbb{S}_{\sigma}(\lambda)$ as a Fiedler pencil of $G(\lam)$.
\end{definition}

Now, we define system linearization of $\mathcal{S}(\lam)$ in (\ref{smh}).

\begin{definition}[System linearization]
A pencil of the form $\mathbb{L}(\lam) = \left[
                               \begin{array}{c|c}
                                 L(\lam) & X \\
                                 \hline
                                 Y & D \\
                               \end{array}
                             \right]
 \in \C^{(nm+r) \times (nm+r)}$, where $L(\lam)$ is an $nm \times nm$ matrix pencil, is called a system linearization of $\mathcal{S}(\lam)$ given in (\ref{smh}) if $\mathcal{S}(\lam)$ is strict system equivalent to $\mathbb{L}(\lam)$.
\end{definition}

Next, our aim is to show that the Fiedler pencils associated with any bijection defined in Definition \ref{dfphos} is a linearization of $\mathcal{S}(\lam). $

\begin{theorem}\label{bfflptflr3h}
Let $\mathcal{S}(\lambda)$ be the system matrix given in (\ref{smh}).
Let $L_{\sigma}(\lambda) = \lambda M_{m} - M_{\sigma}$ be the Fiedler linearization of $P(\lambda)$ associated with a bijection $\sigma :\{0, 1, \ldots, m-1\} \rightarrow \{1, 2, \ldots, m\}$ having CISS$(\sigma) = (c_{1}, i_{1}, \ldots, c_{l}, i_{l})$. Then the Fiedler pencil $\mathbb{S}_\sigma(\lam)$ associated with $\sig$ defined by
$$ \mathbb{S}_{\sigma}(\lambda) := \left[
                             \begin{array}{c|c}
                                -L_{\sigma}(\lambda) & e_m\otimes B \\
                              \hline
                            e_{m -c_{1}}^{T} \otimes C & D \\
                              \end{array}
                                  \right] \indent \mbox{ if } c_1>0 $$ and  $$ \mathbb{S}_{\sigma}(\lambda) := \left[
                             \begin{array}{c|c}
                                -L_{\sigma}(\lambda) & e_{m-i_1}\otimes B \\
                              \hline
                            e_{m}^{T} \otimes C & D \\
                              \end{array}
                                  \right] \indent \mbox{ if } c_1=0 $$ is a system linearization of $\mathcal{S}(\lam)$.
\end{theorem}

\begin{proof}
By Theorem $5.1$ given in \cite{rafinami}, $\mathbb{S}_{\sig}(\lam) = \lam \mathbb{M}_m - \mathbb{M}_{\sig}$ is the Fiedler pencil of $\mathcal{S}(\lam)$ associated with $\sig$.
Now, we have
\begin{equation*}
\left[
  \begin{array}{c|c}
    I_{(m-1)n} & 0 \\
    \hline
    0 & \mathcal{S}(\lam) \\
  \end{array}
\right]
 = \left[
  \begin{array}{cc|c}
    I_{(m-1)n}  &  & 0_{n \times r} \\
     & -P(\lambda) & B \\
    \hline
    0_{r \times n} & C & D \\
  \end{array}
\right].
\end{equation*}
Since
\begin{equation*}
U(\lam)(-L_{\sig}(\lam)) V(\lam) = \left[
                                   \begin{array}{cc}
                                    I_{(m-1)n} &  \\
                                     & -P(\lambda) \\
                                      \end{array}
                                       \right],
\end{equation*}
we have
\begin{align}
&\left[
  \begin{array}{c|c}
  U(\lam)(-L_{\sig}(\lam)) V(\lam) & e_{m}\otimes B \\
\hline
 e_{m}^{T} \otimes C & D \\
  \end{array}
\right]  \nonumber \\
& = \left[
      \begin{array}{c|c}
       U(\lam) & 0 \\
       \hline
        0 & I_{r} \\
      \end{array}
    \right]\left[
             \begin{array}{c|c}
              -L_{\sig}(\lam) & U(\lambda)^{-1}(e_{m}\otimes B) \\
               \hline
           (e_{m}^{T} \otimes C)V(\lambda)^{-1} & D  \\
             \end{array}
           \right]\left[
                    \begin{array}{c|c}
                     V(\lambda) & 0 \\
                      \hline
                     0 & I_{r} \\
                    \end{array}
                  \right]. \label{ser4}
\end{align}
By Theorem $3.11$ and Theorem $3.13$ given in \cite{rafinami2} we have
$$ U(\lam)^{-1} (e_{m} \otimes I_{n}) =
\begin{cases}
e_{m} \otimes I_{n}, & \text{ if } c_{1} >0, \\
e_{m-i_{1}}\otimes I_{n}, & \text{ if } c_{1} = 0,
\end{cases}
$$ and $(e_{m}^{T} \otimes I_{n})V(\lam)^{-1} = e_{m - c_{1}}^{T} \otimes I_{n}$. Hence by (\ref{ser4}) we have
$$\left[
             \begin{array}{c|c}
              -L_{\sig}(\lam) & U(\lambda)^{-1}(e_{m}\otimes B) \\
               \hline
           (e_{m}^{T} \otimes C)V(\lambda)^{-1} & D  \\
             \end{array}
           \right] = \mathbb{S}_{\sig}(\lam) $$
and $$\left[
  \begin{array}{c|c}
    I_{(m-1)n} & 0 \\
    \hline
    0 & \mathcal{S}(\lam) \\
  \end{array}
\right] = \left[
            \begin{array}{c|c}
              U(\lam) &  \\
              \hline
               & I_r \\
            \end{array}
          \right] \mathbb{S}_{\sigma}(\lambda) \left[
                                                 \begin{array}{c|c}
                                                   V(\lam) &  \\
                                                   \hline
                                                    & I_r \\
                                                 \end{array}
                                               \right].
$$ This shows that $\mathbb{S}_{\sigma}(\lambda)$ is a system linearization of $S(\lam)$.
\end{proof}

\begin{exam}\label{comp}
Consider the pencil $\mathcal{C}_{1}(\lam):=\mathbb{S}_{\sig}(\lam) = \lam \mathbb{M}_m -  \mathbb{M}_{m-1} \ldots \mathbb{M}_1 \mathbb{M}_0 = \lam \mathbb{M}_m - \mathbb{M}_{\sig}$ of $\mathcal{S}(\lam)$ given in (\ref{smh}). Then CISS$(\sig) = (0, m-1)$. Since $L_{\sig}(\lam) = \lam M_m - M_{m-1} \cdots M_{1} M_{0}$ is the first companion form of $P(\lam)$, by the Theorem \ref{bfflptflr3h} we have
$$
\mathbb{S}_{\sig}(\lam) =  \left[
                          \begin{array}{c|c}
                            L_{\sig}(\lam) & e_1 \otimes B \\
                            \hline
                            e_m^{T} \otimes C & D \\
                          \end{array}
                        \right] $$
               $$ = \lam \left[
                         \begin{array}{cccc|c}
                           A_{m} &  &  &  &  \\
                            & I_{n} &  &  &  \\
                            &  & \ddots &  &  \\
                            &  &  & I_{n} &  \\
                            \hline
                            &  &  &  & 0_r \\
                         \end{array}
                       \right]- \left[
                \begin{array}{cccc|c}
                  -A_{m-1} & -A_{m-2}& \cdots & -A_{0}& -B  \\
                   I_{n} & 0 & \cdots & 0 &  \\
                   & \ddots &  & \vdots  & \\
                   &  & I_{n} & 0 & \\
                   \hline
                  &  &  & -C &  -D \\
                \end{array}
              \right]
$$
is a system linearization of $\mathcal{S}(\lam)$. Therefore we refer to $\mathcal{C}_{1}(\lam)$ as {\em first companion form} of $\mathcal{S}(\lam)$.

Consider $\mathcal{C}_{2}(\lam) :=\mathbb{S}_{\sig}(\lambda) = \lambda \mathbb{M}_{m} - \mathbb{M}_{0}\mathbb{M}_{1} \cdots \mathbb{M}_{m-2}\mathbb{M}_{m-1} =\lam \mathbb{M}_{m} - \mathbb{M}_{rev\,\sig}. $ Then CISS$(\sig) = (m-1, 0)$. Since $L_{\sig}(\lam) = \lam M_m - M_{0}M_{1}\cdots M_{m-1}$ is the second companion form of $P(\lam)$, by the Theorem \ref{bfflptflr3h} we have
$$\mathcal{C}_{2}(\lam) = \lam \left[
                                    \begin{array}{cccc|c}
                                      A_{m} &  &  &  &  \\
                                       & I_{n} &  &  &  \\
                                       &  & \ddots &  &  \\
                                       &  &  & I_{n} &  \\
                                       \hline
                                       &  &  &  & 0 \\
                                    \end{array}
                                  \right] - \left[
                                                 \begin{array}{cccc|c}
                                                 -A_{m-1} & I_{n} &  &  &   \\
                                                 -A_{m-2} & 0 & \ddots &  &  \\
                                                 \vdots & \ddots &  & I_{n} & \\
                                                 -A_{0} & \cdots & 0 & 0  & -B \\
                                                 \hline
                                                 -C  &    &    &   &  -D \\
                                           \end{array}
                                           \right].$$
So, we refer  $\mathcal{C}_{2}(\lambda)$ as second companion form of $\mathcal{S}(\lam)$ or $G(\lam)$.
$\blacksquare$
\end{exam}

It is easy to see that if $\lam$ is an eigenvalue of $G(\lam)$ then $G(\lam)x = 0$ if and only if $\mathcal{C}_{1}(\lambda)v = 0, $
where \begin{equation}
v= \left[
     \begin{array}{c}
       \lambda^{m-1}x \\
       \lambda^{m-2}x \\
       \vdots \\
       x \\
       \hline
       (P(\lam))^{-1}Bx \\
        \end{array}
        \right]. \label{revofcfotf}
\end{equation}
Equivalently, let $\mathcal{N}(G(\lam))$ and $\mathcal{N}(\mathcal{C}_1(\lam))$ be nullspaces of $G(\lam)$ and $\mathcal{C}_1(\lam)$ respectively, then the map
$$
\mathcal{N}(G(\lam))  \,\, \rightarrow \mathcal{N}(\mathcal{C}_1(\lam)), \indent
 x \mapsto \left[
     \begin{array}{c}
       \lambda^{m-1}x \\
       \lambda^{m-2}x \\
       \vdots \\
       x \\
       \hline
       (P(\lam))^{-1}Bx \\
        \end{array}
        \right]
$$ is an isomorphism.

Now, we construct low band-width linearization of $\mathcal{S}(\lam). $

\begin{theorem}\label{pdpchos}
Let $L_{\sigma}(\lambda) = \lambda M_{m} - M_{\sigma}$ be a pentadiagonal Fiedler pencil of $P(\lambda)$ associated with a bijection $\sigma :\{0, 1, \ldots, m-1\} \rightarrow \{1, 2, \ldots, m\}$ having CISS$(\sigma) = (c_{1}, i_{1}, \ldots, c_{l}, i_{l})$. Then $\mathbb{S}_\sigma(\lam)= \lambda \mathbb{M}_{m} - \mathbb{M}_{\sigma}$ is pentadiagonal Fiedler pencil of $\mathcal{S}(\lam) $ provided that $c_1 \leq 1$ and $i_1 \leq 1$. In such a case, we have
$$ \mathbb{S}_{\sigma}(\lambda) := \left[
                             \begin{array}{c|c}
                                -L_{\sigma}(\lambda) & e_m\otimes B \\
                              \hline
                            e_{m -1}^{T} \otimes C & D \\
                              \end{array}
                                  \right] \indent \mbox{ if } c_1>0 $$ $$\mbox{ and }\indent  \mathbb{S}_{\sigma}(\lambda) := \left[
                             \begin{array}{c|c}
                                -L_{\sigma}(\lambda) & e_{m-1}\otimes B \\
                              \hline
                            e_{m}^{T} \otimes C & D \\
                              \end{array}
                                  \right] \indent \mbox{ if } c_1=0. $$
\end{theorem}

\begin{proof}
If $c_1>0$ then $c_1 = 1$ and $ e_{m -c_{1}}^{T} \otimes C = e_{m -1}^{T} \otimes C$. So by the Theorem \ref{bfflptflr3h} we have  $$\mathbb{S}_{\sigma}(\lambda) =\left[
                             \begin{array}{c|c}
                                -L_{\sigma}(\lambda) & e_m\otimes B \\
                              \hline
                            e_{m -c_{1}}^{T} \otimes C & D \\
                              \end{array}
                                  \right]
                                    = \left[
                             \begin{array}{c|c}
                                -L_{\sigma}(\lambda) & e_m\otimes B \\
                              \hline
                            e_{m -1}^{T} \otimes C & D \\
                              \end{array}
                                  \right] $$ is a pentadiagonal,
since $L_{\sigma}(\lambda)$ is pentadiagonal. Similar argument holds for $c_1 = 0$ and $i_1 = 1. $
\end{proof}

\begin{exam}
Let $P(\lam) = \lam^6 A_6 + \cdots + \lam A_1 +A_0$ be regular and $G(\lam) = C P(\lam)^{-1}B + D$. Consider the pencil $\mathbb{S}_{\sig}(\lam) = \lam \mathbb{M}_6 - \mathbb{M}_{1}\mathbb{M}_{3}\mathbb{M}_{5}\mathbb{M}_{0}\mathbb{M}_{2}\mathbb{M}_{4} = \lambda \mathbb{M}_{6} - \mathbb{M}_{\sig}. $ Now
$$\mathbb{M}_{\sigma} = \left[
    \begin{array}{cccccc|c}
      -A_{5} & -A_{4} & I_{n} & 0 & 0 & 0 & 0 \\
      I_{n} & 0 & 0 & 0 & 0 & 0 & 0 \\
      0 & -A_{3} & 0 & -A_{2} & I_{n} & 0 & 0 \\
      0 & I_{n} & 0 & 0 & 0 & 0 & 0 \\
      0 & 0 & 0 & -A_{1} & 0 & -A_{0} & -C \\
      0 & 0 & 0 & I_{n} & 0 & 0 & 0 \\
       \hline
      0 & 0 & 0 & 0 & 0 & -B & -D \\
    \end{array}
  \right]
$$
shows that $\mathbb{S}_{\sig}(\lam)$ is pentadiagonal. $\blacksquare$
\end{exam}

Next, we consider the higher order LTI state space system and associated system matrix. We
derive a linearized state space system of the given higher order system and show that both the systems are strict system equivalent with each other.

Consider the system $\Sigma_1$ given in (\ref{mulsss}). Then
the system $\Sigma$ can be written as
$$\mathcal{S}(\lam) \left[
                      \begin{array}{c}
                        x \\
                        u \\
                      \end{array}
                    \right] = \left[
                                \begin{array}{c}
                                  0 \\
                                  y \\
                                \end{array}
                              \right].
$$
The nonzero vector $u$ for which $G(\lam)u = 0$ is called a {\em zero direction} of $\Sigma_1$. The nonzero vector $\left[
                                   \begin{array}{c}
                                     x \\
                                     u \\
                                   \end{array}
                                 \right]$ such that $\mathcal{S}(\lam) \left[
                                                                         \begin{array}{c}
                                                                           x \\
                                                                           u \\
                                                                         \end{array}
                                                                       \right] = 0 $
is called an {\em invariant direction} of $\Sigma_1,$ where $u$ is referred as {\em input zero direction} and $x$ is referred to as {\em state zero direction} associated with $\lam$ \cite{vardulakis}. Now, we derive the linearized state-space system.

\begin{theorem}\label{lsssohos}
Let $\Sigma_1$ be the LTI state-space system given in (\ref{mulsss}) and $\mathcal{S}(\lam)$ be the associated system matrix in (\ref{smh}). Let $E_{\sig}(P), F_{\sig}(P)$ and $ K_{\sig}(P)$ be as in Theorem \ref{ralitfmp} and Theorem \ref{dsigp} respectively. If $\mathbb{S}_{\sig}(\lam)$ is the Fiedler pencil of $\mathcal{S}(\lam)$ associated with a bijection $\sig$, then the state-space system $\Sigma_2$ associated with $\mathbb{S}_{\sig}(\lam)$ is given by
\begin{eqnarray*}
\lam M_m X &=& M_{\sig} X +  K_{\sig}(P)^{T}B u  \\
y &=& C F_{\sig}(P) X + Du
\end{eqnarray*}
where $L_{\sig}(\lam) = \lam M_m - M_{\sig}$ is the Fiedler pencil of $P(\lam)$ associated with the bijection $\sig$.
The state vector $x$ of $\Sigma_1$ is given by $x= F_{\sig}(P)X$.
Further, if CISS $(\sig) = (c_1, i_1, \ldots, c_l, i_l)$, then $\Sigma_2$ is given by
\begin{equation*}
\begin{aligned}
\lam M_m X &= M_{\sig} X + (e_m \otimes B) u, \\
y &= (e_{m-c_1}^{T} \otimes C) X + D u
\end{aligned} \indent \mbox{ if } c_1 > 0
\end{equation*}
and
\begin{equation*}
\begin{array}{c}
\lam M_m X = M_{\sig} X+ (e_{m-i_1} \otimes B) u, \\
y = (e_{m}^{T} \otimes C) X + D u
\end{array} \indent \mbox{ if } c_1 = 0.
\end{equation*}
Further, the state vector $x$ of $\Sigma_1$ is given by
$x = (e_{(m-c_1)}^{T} \otimes I_n) X$. Furthermore, the system $\Sigma_1$ is strict system equivalent to the system $\Sigma_2. $
\end{theorem}

\begin{proof}
We have
\begin{align*}
\left[
  \begin{array}{c|c}
    -P(\lam) & B \\
    \hline
    C & D \\
  \end{array}
\right] \left[
          \begin{array}{c}
            x \\
            u \\
          \end{array}
        \right] = \left[
                    \begin{array}{c}
                      0 \\
                      y \\
                    \end{array}
                  \right].
\end{align*}
Considering the extended system, we have
$$
 \left[
  \begin{array}{c|cc}
    I_{(m-1)n} &  &  \\
    \hline
     & -P(\lam) & B \\
     & C & D \\
  \end{array}
\right] \left[
          \begin{array}{c}
            \left[
              \begin{array}{c}
                0 \\
                x \\
              \end{array}
            \right]\\
            u \\
          \end{array}
        \right] = \left[
                    \begin{array}{c}
                      \left[
                        \begin{array}{c}
                          0 \\
                          0 \\
                        \end{array}
                      \right]\\
                      y \\
                    \end{array}
                  \right] $$
$$ = \left[
      \begin{array}{c|c}
        U(\lam) &  \\
        \hline
         & I_r \\
      \end{array}
    \right]\mathbb{S}_{\sig}(\lam) \left[
                            \begin{array}{c|c}
                              V(\lam) &  \\
                              \hline
                               & I_r \\
                            \end{array}
                          \right] \left[
                                    \begin{array}{c}
                                      \left[
                                        \begin{array}{c}
                                          0 \\
                                          x \\
                                        \end{array}
                                      \right]\\
                                      u \\
                                    \end{array}
                                  \right] = \left[
                                              \begin{array}{c}
                                                \left[
                                                  \begin{array}{c}
                                                    0 \\
                                                    0 \\
                                                  \end{array}
                                                \right] \\
                                                y \\
                                              \end{array}
                                            \right]$$
$$ \Rightarrow \mathbb{S}_{\sig}(\lam)\left[
                                       \begin{array}{c}
                                         V(\lam)(e_m \otimes I_n)x \\
                                         u \\
                                       \end{array}
                                     \right] = \left[
                                                 \begin{array}{c|c}
                                                   U(\lam)^{-1} & 0 \\
                                                   \hline
                                                   0 & I_r \\
                                                 \end{array}
                                               \right] \left[
                                                         \begin{array}{c}
                                                           0 \\
                                                           y \\
                                                         \end{array}
                                                       \right] $$
$$ \Rightarrow \mathbb{S}_{\sig}(\lam)\left[
                                       \begin{array}{c}
                                         E_{\sig}(P)x \\
                                         u \\
                                       \end{array}
                                     \right] = \left[
                                                 \begin{array}{c}
                                                   0 \\
                                                   y \\
                                                 \end{array}
                                               \right] $$
$$ \Rightarrow \left[
              \begin{array}{c|c}
                -L_{\sig}(\lam) & K_{\sig}(P)^{T}B \\
                \hline
                C F_{\sig}(P) & D \\
              \end{array}
            \right] \left[
                                       \begin{array}{c}
                                         E_{\sig}(P)x \\
                                         u \\
                                       \end{array}
                                     \right] = \left[
                                                 \begin{array}{c}
                                                   0 \\
                                                   y \\
                                                 \end{array}
                                               \right].
$$
This gives
\begin{eqnarray*}
 -L_{\sig}(\lam)(E_{\sig}(P)x) + K_{\sig}(P)^{T}B u &=& 0  \\
 C F_{\sig}(P) (E_{\sig}(P)x) + Du &=& y.
\end{eqnarray*}
Thus the state-space system associated with $\mathbb{S}_{\sig}(\lam)$ is given by
\begin{eqnarray*}
-L_{\sig}(\lam)X + K_{\sig}(P)^{T}B u &=& 0  \\
 C F_{\sig}(P) X + Du &=& y.
\end{eqnarray*}
where $X :=  E_{\sig}(P)x$. Hence
\begin{equation*}
\begin{aligned}
\lam M_m X &= M_{\sig} X +  K_{\sig}(P)^{T}B u  \\
y &= C F_{\sig}(P) X + Du
\end{aligned}
\end{equation*}
is the required state-space system associated with $\mathbb{S}_{\sig}(\lam). $ Since $F_{\sig}(P) E_{\sig}(P) = I_n$, we have $F_{\sig}(P) X = F_{\sig}(P) E_{\sig}(P)x = x$. Second part of the proof follows from Theorem \ref{bfflptflr3h}.
\end{proof}

Consider the system $\Sigma_2$ given in Theorem \ref{lsssohos}. Then the associated transfer function $\mathbb{G}(\lam)$ is given by
$$\mathbb{G}(\lam) = (e_{m-c_1}^{T} \otimes C) L_{\sig}(\lam)^{-1} (e_m \otimes B) + D,  \indent \mbox{ if } c_1>0$$
and
$$\mathbb{G}(\lam) = (e_{m}^{T} \otimes C) L_{\sig}(\lam)^{-1} (e_{m-i_1} \otimes B) + D  \indent \mbox{ if } c_1=0. $$
It is easy to see that $\mathbb{G}(\lam) = G(\lam)$, which is consistent with the fact that $\mathcal{S}(\lam)$ is SSE to $\mathbb{S}_{\sig}(\lam)$.

\begin{theorem}
Let $\Sigma_1$ and $\Sigma_2$ be as in Theorem \ref{lsssohos}. Let $G(\lam)$ and $\mathbb{G}(\lam)$ be associated transfer functions of  $\Sigma_1$ and $\Sigma_2,$ respectively. Then we have the following.
\begin{itemize}

\item [(a)] The system $\Sigma_1$ is observable if and only if the system $\Sigma_2$ is observable.

\item[(b)] The system $\Sigma_1$ is controllable if and only if the system $\Sigma_2$ is controllable.

\item [(c)] The transfer function $\mathbb{G}(\lam)$  is unimodularly equivalent to $\left[
               \begin{array}{cc}
                 I_{(m-1)n} &  \\
                  & G(\lam) \\
               \end{array}
             \right]$.

\item [(d)] Let $\lam \in \C$. Then $\lam$ is an input (output, input-output)  decoupling zero of $\Sigma_1$ if and only if $\lam$ is an input (output, input-output)  decoupling zero of $\Sigma_2$.

\end{itemize}
\end{theorem}

\begin{proof}
Since $\rank[L_{\sig}(\lam),\, B] = \rank[L_{\sig}(\lam), \, B F_{\sig}(P)]$ and $$\rank\left[
                                                                          \begin{array}{c}
                                                                            L_{\sig}(\lam) \\
                                                                            C \\
                                                                          \end{array}
                                                                        \right] = \rank \left[
                                                                                          \begin{array}{c}
                                                                                            L_{\sig}(\lam)\\
                                                                                            K_{\sig}(P)^{T} C \\
                                                                                          \end{array}
                                                                                        \right]
$$ for all $\lam \in \C$, the desired results in $(a), (b)$ and $(d)$ follow. The result in $(c)$ is immediate.
\end{proof}

Next, we define proper generalized Fiedler (PGF) pencils of $G(\lam)$ given in (\ref{tfh}). These pencils are defined as follows.

\begin{definition}[PGF pencil] A permutation $ w :=( w_0, w_1) $ of $\{0, 1, \ldots, m\}$ is said to be proper if $ 0 \in w_0$ and $ m \in w_1.$ Let $ w := (w_0, w_1)$ be a proper permutation of  $\{0, 1, \ldots, m\}$  and define $\widehat{\mathbb{M}}_{w_0} := \mathbb{M}_{w_0} $ and
$$ \widehat{\mathbb{M}}_{w_1} := \mathbb{M}_{i_1}^{-1} \cdots \mathbb{M}_{i_p}^{-1} \mathbb{M}_m \mathbb{M}_{j_1}^{-1} \cdots \mathbb{M}_{j_q}^{-1},  \,\, \mbox{ if } w_1 =( i_1, \cdots, i_p, m, j_1, \cdots, j_q).$$
 Then the pencil $\mathbb{K}_w(\lam) := \lam \widehat{\mathbb{M}}_{w_1} - \widehat{\mathbb{M}}_{w_0}$ is called a proper generalized Fiedler (PGF) pencil of $G(\lam)$ associated with $w.$  We also refer to $\mathbb{K}_w(\lam)$ as the PGF pencil of $\mathcal{S}(\lam)$ given in (\ref{smh}) associated with $w.$
\end{definition}

Clearly a Fiedler pencil $\mathbb{S}_\sigma(\lam)$ is a PGF pencil of $G(\lam).$ Indeed, for $ w :=(\sigma, m)$ we have $ \mathbb{K}_w (\lam)= \mathbb{S}_{\sigma}(\lam).$

A product $\mathbb{M}_{\textbf{q}}$ corresponding to the index tuple $\textbf{q} = (i_1, i_2, \ldots, i_m)$ is said to be {\em operation-free} if the block entries (up to sign) of $\mathbb{M}_{\textbf{q}}$ consist of matrices from  $0, I_{n}, I_{r}, C, B, A, E$ and  $A_0, A_1, \ldots, A_m$, \cite{rafinami3}. Thus, it is easily seen that a PGF pencil $ \mathbb{K}_w(\lam) = \lam \widehat{\mathbb{M}}_{w_1} - \widehat{\mathbb{M}}_{w_0}$ of $G(\lam)$ associated with a proper permutation $w = (w_0, w_1)$ of $\{0, 1, \ldots, m\}$ is operation-free.

\begin{remark}
Note that $\mathbb{M}_m$ is always singular. Also note that non-proper generalized Fiedler (NPGF) pencils of $G(\lam)$ are defined only when  $\mathbb{M}_{m}$ and/or $\mathbb{M}_{0}$ are nonsingular. Thus unlike NPGF pencils for LTI state space system in state space form, here it is impossible to define the NPGF pencils for higher order systems.
\end{remark}

Next, we show that  a  PGF pencil of $G(\lam)$ given in (\ref{tfh}) can be constructed from a PGF pencil $P(\lam)$ and vice-versa. In fact, there is a bijection from the set of PGF pencils of $P(\lam)$ to the set of PGF pencils of $G(\lam).$

\begin{theorem}[PGF Pencil]\label{pgf} Let $w := (w_0, w_1)$ be a proper permutation of $\{0, 1, \ldots, m\}.$ Let $\mathbb{K}_{w}(\lam) := \lam \widehat{\mathbb{M}}_{w_1}-\widehat{\mathbb{M}}_{w_0} $ and $K_{w}(\lam) := \lam \widehat{M}_{w_1} -  \widehat{M}_{w_0}$  be the  PGF pencils of $G(\lam)$ and $P(\lam)$ associated with $w,$ respectively.  If CIP$(w_0) =(c_{0}, i_{0})$ then
$$  \mathbb{K}_{w}(\lam)  =
 \left[\begin{array}{c|c} K_{w}(\lam) & e_{m-i_0}\otimes B \\ \hline
 e_{m-c_0}^T \otimes C & D \\ \end{array} \right].$$
Thus,  the map $$\mathrm{PGF}(P) \longrightarrow \mathrm{PGF}(G), \, K_w(\lam) \longmapsto  \left[\begin{array}{c|c} K_{w}(\lam) & e_{m-i_0}\otimes B \\ \hline
 e_{m-c_0}^T \otimes C & D \\ \end{array} \right] $$
 is a bijection, where $\mathrm{PGF}(P)$ and $\mathrm{PGF}(G)$ denote the set of PGF pencils of $P(\lam)$ and $G(\lam)$, respectively.
\end{theorem}

\begin{proof}
The proof is directly follows from proof of the Theorem $2.7$ in \cite{rafinami3}.
\end{proof}

The next result characterizes the block tridiagonal PGF pencils of $G(\lam)$ given in (\ref{tfh}).

\begin{theorem}\label{tridpgf}  Let  $K_{w}(\lambda) := \lambda \widehat{M}_{w_1} - \widehat{M}_{w_0}$ and $\mathbb{K}_{w}(\lam) := \lam \widehat{\mathbb{M}}_{w_1} - \widehat{\mathbb{M}}_{w_0},$ respectively, be the PGF pencils of $P(\lam)$ and $G(\lam)$ associated with a proper permutation $w:=(w_0, w_1)$ of $ \{0, 1, \ldots, m\}.$
Then $ \mathbb{K}_{w}(\lam)$ is block tridiagonal  if and only if $K_{w}(\lam)$ is block tridiagonal and CIP$(w_0) = (0,0). $ In such a case, the block tridiagonal pencil is given by
\be\label{tridiag}  \mathbb{K}_w(\lambda) = \left[
                             \begin{array}{c|c}
                                K_{w}(\lambda) & e_m\otimes B \\
                              \hline
                            e_{m}^{T} \otimes C & D \\
                              \end{array}
                                  \right]. \ee

\end{theorem}

\begin{proof}
The proof is directly follows from proof of the Theorem $2.10$ in \cite{rafinami3}.
\end{proof}

We now construct a block tridiagonal PGF pencils of $G(\lam)$  which may be useful for computation.

\begin{theorem}\label{tridiagp}
Consider the sub-permutations  $\sig_{1} := (1, 3, 5, \ldots)$ and $\sig_{2} := (2, 4, 6, \ldots)$ of $\{ 1, 2, \ldots, m-1\}.$
Then the PGF pencil $ \mathbb{K}_o(\lam) := \lam \mathbb{M}_{\sig_{1}}^{-1}\mathbb{M}_{m} - \mathbb{M}_{0}\mathbb{M}_{\sig_{2}}$ of $G(\lam)$ is  block tridiagonal. The PGF pencil $ \mathbb{K}_e (\lam) := \lam \mathbb{M}_{m} \mathbb{M}_{\sig_{2}}^{-1}- \mathbb{M}_{\sig_{1}} \mathbb{M}_{0}$ is not block tridiagonal.
\end{theorem}

\begin{proof}
The proof is directly follows from proof of the Theorem $2.11$ in \cite{rafinami3}.
\end{proof}

\subsection{Hermitian PGF pencil}
Consider $\mathcal{S}(\lam)$ and $G(\lam)$ are given in (\ref{smh}) and (\ref{tfh}), respectively. We define the adjoint of $\mathcal{S}(\lam)$  by
$$\mathcal{S}^{*}(\lam) = \left[
                                         \begin{array}{c|c}
                                           -P^{*}(\lam) & C^{*} \\
                                           \hline
                                           B^{*} & D^{*} \\
                                         \end{array}
                                       \right],
$$ where $P^{*}(\lam) = \sum_{j = 1}^{m}\lam^{j}A_j^{*}$. The transposes of $\mathcal{S}(\lam)$ and $P(\lam)$ are defined similarly. The system matrix $\mathcal{S}(\lam)$ is said to be Hermitian if  $\mathcal{S}^*(\lam) = \mathcal{S}(\lam)$ for $\lam \in \C$. The transfer function  $G(\lam)$  is said to be Hermitian (resp., symmetric) when $\mathcal{S}(\lam)$ is
Hermitian (resp., symmetric).

\begin{theorem}\label{sym}  Let  $K_{w}(\lambda) := \lambda \widehat{M}_{w_1} - \widehat{M}_{w_0}$ and $\mathbb{K}_{w}(\lam) := \lam \widehat{\mathbb{M}}_{w_1} - \widehat{\mathbb{M}}_{w_0}$  be the PGF pencils of $P(\lam)$ and $G(\lam),$ respectively, associated with a proper permutation $w:=(w_0, w_1)$ of $ \{0, 1, \ldots, m\}.$ Suppose that $G(\lam)$ is Hermitian.  Then $ \mathbb{K}_{w}(\lam)$ is  Hermitian  if and only if $K_{w}(\lam)$ is  Hermitian and CIP$(w_0) = (0,0). $ In such a case, the Hermitian pencil $\mathbb{K}_{w}(\lam)$  is given by  (\ref{tridiag}) with $C = B^*.$
\end{theorem}
\begin{proof}
By  Theorem~\ref{pgf},  $ \mathbb{K}_{w}(\lam)$ is Hermitian if and only if $K_w(\lam)$ is Hermitian and $(c_0, i_0) = (0,0).$ Thus
$ \mathbb{K}_{w}(\lam)$ is given by $(\ref{tridiag})$ with $K_w(\lam)$ being Hermitian and $ C= B^*.$
\end{proof}

We now show that a  Hermitian $G(\lam)$  admits a  Hermitian GF pencil when the degree $m$ of the  matrix polynomial $P(\lam)$  is an odd number.

\begin{theorem}\label{slormf}
Let $G(\lam)$ be a Hermitian (resp., symmetric). Suppose that  $m,$ the degree of $P(\lam)$, is odd. Then the PGF pencil $\mathbb{K}(\lam)$ of $G(\lam)$  given by
$$ \mathbb{K}(\lambda) =
\lambda \mathbb{M}_{m}\mathbb{M}_{m-2}^{-1}\cdots \mathbb{M}_{3}^{-1}\mathbb{M}_{1}^{-1} - \mathbb{M}_{0}\mathbb{M}_{2} \cdots \mathbb{M}_{m-3}\mathbb{M}_{m-1} \label{spr}
$$
is Hermitian (resp., symmetric). Thus  $\mathbb{K}(\lam)$ is a Hermitian (resp., symmetric)  linearization of $\mathcal{S}(\lam)$.
\end{theorem}

\begin{proof}
The proof is directly follows from proof of the Theorem~4.2 in \cite{rafinami3}.
\end{proof}

\subsection{ Zero direction recovery from linearizations} Our main aim in this section is that given an eigenvector $x$ of $\mathbb{K}(\lam)$ determine an eigenvector of $\mathcal{S}(\lam)$ from $x$. That is, to recover eigenvectors of $G(\lam)$ and $\mathcal{S}(\lam)$ from those of the Fiedler pencils and the PGF pencils of $G(\lam).$ For the rest of this section, we assume that $G(\lam)$ is regular. Hence $\mathcal{S}(\lam)$ is also regular. As $G(\lam)$ is regular, recall that $ \lam$ is an eigenvalue of $G(\lam)$ provided that $\rank(G(\lam)) < n$ and that  $\mathcal{N}_r(G(\lam))$ and $\mathcal{N}_l(G(\lam))$, respectively, are the right and left null spaces of $G(\lam).$ Now, we discuss recovery of zero directions of $\Sigma_1$ from eigenvector recovery of $\mathbb{S}_{\sig}(\lam)$.

\begin{theorem} \cite{behera, rafinami2} \label{evrthfmpdm}
Let $\mathbb{L}_{\sig}(\lam)$ be the Fiedler linearization of $\mathcal{S}(\lam)$ associated with a bijection $\sig$. Let $\lam \in \C$ and $E_{\sig}(P)$, and $H_{\sig}(P)$ be as in Theorem~\ref{ralitfmp} and \ref{dsigp}, respectively. Define $\mathbb{E}_{\sig}(\mathcal{S}) : \C^{n+r}\rightarrow \C^{nm+r}$ and  $\mathbb{H}_{\sig}(\mathcal{S}) : \C^{n+r}\rightarrow \C^{nm+r}$ by
$$\mathbb{E}_{\sig}(\mathcal{S})= \left[
                                    \begin{array}{c|c}
                                      E_{\sig}(P) & 0 \\
                                      \hline
                                      0 & I_r \\
                                    \end{array}
                                  \right]
 \text{ and } \mathbb{H}_{\sig}(\mathcal{S})= \left[
                                    \begin{array}{c|c}
                                      H_{\sig}(P) & 0 \\
                                      \hline
                                      0 & I_r \\
                                    \end{array}
                                  \right]. $$ Then $\mathbb{E}_{\sig}(\mathcal{S}) : \mathcal{N}_r(\mathcal{S}(\lam)) \rightarrow \mathcal{N}_r(\mathbb{L}_{\sig}(\lam))$ and $\mathbb{H}_{\sig}(\mathcal{S}) : \mathcal{N}_l(\mathcal{S}(\lam)) \rightarrow \mathcal{N}_l(\mathbb{L}_{\sig}(\lam))$ are isomorphisms.

Let $F_{\sig}(P)$, and $K_{\sig}(P)$ be as in Theorem~\ref{ralitfmp} and \ref{dsigp}, respectively. Define $\mathbb{F}_{\sig}(\mathcal{S}) : \C^{n+r}\rightarrow \C^{nm+r}$ and  $\mathbb{K}_{\sig}(\mathcal{S}) : \C^{n+r}\rightarrow \C^{nm+r}$ by
$$\mathbb{F}_{\sig}(\mathcal{S})= \left[
\begin{array}{c|c}
F_{\sig}(P) & 0 \\
\hline
0 & I_r \\
\end{array}
\right]
\text{ and } \mathbb{K}_{\sig}(\mathcal{S})= \left[
\begin{array}{c|c}
K_{\sig}(P) & 0 \\
\hline
0 & I_r \\
\end{array}
\right]. $$ Then $\mathbb{F}_{\sig}(\mathcal{S}) : \mathcal{N}_r(\mathbb{L}_{\sig}(\lam)) \rightarrow \mathcal{N}_r(\mathcal{S}(\lam))$ and $\mathbb{K}_{\sig}(\mathcal{S}) : \mathcal{N}_l(\mathbb{L}_{\sig}(\lam)) \rightarrow \mathcal{N}_l(\mathcal{S}(\lam))$ are isomorphisms.
\end{theorem}

Next, we consider Fiedler linearization $\mathbb{S}_{\sig}(\lam)$ of the system matrix $\mathcal{S}(\lam)$ in (\ref{smh}) and determine isomorphism between nullspaces of $\mathcal{S}(\lam)$ and $\mathbb{S}_{\sig}(\lam). $
Now, by Theorem~\ref{evrthfmpdm}, it is also easy to see that
$\mathbb{E}_{\sig}(\mathcal{S}) : \mathcal{N}_r(\mathcal{S}(\lam)) \rightarrow \mathcal{N}_r(\mathbb{S}_{\sig}(\lam))$ and $\mathbb{F}_{\sig}(\mathcal{S}) : \mathcal{N}_r(\mathbb{S}_{\sig}(\lam)) \rightarrow \mathcal{N}_r(\mathcal{S}(\lam))$ are isomorphisms.
Furthermore, $\mathbb{H}_{\sig}(\mathcal{S}) : \mathcal{N}_l(\mathcal{S}(\lam)) \rightarrow \mathcal{N}_l(\mathbb{S}_{\sig}(\lam))$ and $\mathbb{K}_{\sig}(\mathcal{S}) : \mathcal{N}_l(\mathbb{S}_{\sig}(\lam)) \rightarrow \mathcal{N}_l(\mathcal{S}(\lam))$ are isomorphisms.
The next result gives an isomorphism between nullspaces of  $G(\lam)$ and $\mathcal{S}(\lam)$.

\begin{theorem}\label{isoogtosh}
Let $\mathcal{S}(\lam)$ and $G(\lam)$ be as in (\ref{smh}) and (\ref{tfh}). Let $\lam \in \C$ be an eigenvalue of $G(\lam)$. Define $f: \C^{r} \rightarrow \C^{n+r}$ by $f(x) = \left[
                                                         \begin{array}{c}
                                                           P(\lam)^{-1}Bx \\
                                                            x \\
                                                         \end{array}
                                                       \right]$ and $g : \C^{r} \rightarrow \C^{n+r}$ by
$g(x) = \left[
           \begin{array}{c}
             C P(\lam)^{-1}x \\
             x \\
           \end{array}
           \right]$. Then the maps
$f \,:  \,\, \mathcal{N}_{r}(G(\lam))\rightarrow \mathcal{N}_{r}(\mathcal{S}(\lam))$  and $g \, :  \,\, \mathcal{N}_l(G(\lam))\rightarrow \mathcal{N}_l(\mathcal{S}(\lam))$ are isomorphisms.
\end{theorem}

\begin{proof}
Let $x \in \mathcal{N}_{r}(G(\lam)) \Rightarrow G(\lam)x = 0$. Now \begin{align*}\mathcal{S}(\lam) \left[
                                                         \begin{array}{c}
                                                            P(\lam)^{-1}Bx \\
                                                            x \\
                                                         \end{array}
                                                       \right] &= \left[
                                                                   \begin{array}{c|c}
                                                                     -P(\lam) & B \\
                                                                     \hline
                                                                     C & D \\
                                                                   \end{array}
                                                                 \right] \left[
                                                         \begin{array}{c}
                                                            P(\lam)^{-1}Bx \\
                                                            x \\
                                                         \end{array}
                                                       \right] \\ &= \left[
                                                                   \begin{array}{c}
                                                                     -Bx+Bx  \\
                                                                     G(\lam)x \\
                                                                   \end{array}
                                                                 \right] = \left[
                                                                             \begin{array}{c}
                                                                               0 \\
                                                                               0 \\
                                                                             \end{array}
                                                                           \right].
 \end{align*}                                       So $\left[
                                                         \begin{array}{c}
                                                           P(\lam)^{-1}Bx \\
                                                           x \\
                                                         \end{array}
                                                       \right] \in \mathcal{N}_{r}(\mathcal{S}(\lam)).$
This shows that  $f: \mathcal{N}_{r}(G(\lam))\rightarrow \mathcal{N}_{r}(\mathcal{S}(\lam)).$ It is easy to check $f: \mathcal{N}_{r}(G(\lam))\rightarrow \mathcal{N}_{r}(\mathcal{S}(\lam))$ is an isomorphism. The proof is similar for $g.$
\end{proof}

Consider $G(\lam)$ and $S(\lam)$ given in (\ref{tfh}) and (\ref{smh}) respectively. Then we show that there is an isomorphism between nullspaces of $G(\lam)$ and a Fiedler linearization of $\mathcal{S}(\lam)$.

\begin{theorem}\label{cevlthos}
Let $\mathbb{S}_{\sig}(\lam)$ be the linearization of $\mathcal{S}(\lam)$ associated with a bijection $\sig$. Let $E_{\sig}(P)$ be as in Theorem \ref{ralitfmp} and $H_{\sig}(P)$ be as in Theorem \ref{dsigp}. Suppose that $G(\lam)$ is given in (\ref{tfh}) and $\lam \in \C$ is an eigenvalue of $G(\lam)$.
\begin{itemize}

\item [(a)] Then $\mathbb{E}_{\sig}(G) : \C^{r} \rightarrow \C^{nm +r} $ defined by
$$\mathbb{E}_{\sig}(G)x = \left[
                                  \begin{array}{c}
                                    E_{\sig}(P) P(\lam)^{-1}Bx \\
                                    x \\
                                  \end{array}
                                \right]
$$
is an isomorphism from $ \mathcal{N}_r(G(\lam))$ to $\mathcal{N}_r(\mathbb{S}_{\sig}(\lam))$.

\item [(b)] Define $\mathbb{H}_{\sig}(G) : \C^{r} \rightarrow \C^{nm +r} $ by
$$\mathbb{H}_{\sig}(G)y = \left[
                                  \begin{array}{c}
                                    H_{\sig}(P) C P(\lam)^{-1}y \\
                                    y \\
                                  \end{array}
                                \right].
$$
Then $\mathbb{H}_{\sig}(G)$ is an isomorphism from $\mathcal{N}_l(G(\lam))$ to $\mathcal{N}_l(\mathbb{S}_{\sig}(\lam)). $
\end{itemize}
\end{theorem}

\begin{proof} By Theorem \ref{isoogtosh}, we have $f: \mathcal{N}_{r}(G(\lam)) \rightarrow \mathcal{N}_{r}(\mathcal{S}(\lam)) $ and $g : \mathcal{N}_{l}(G(\lam)) \rightarrow \mathcal{N}_{l}(\mathcal{S}(\lam))$ are isomorphisms. Since $\mathbb{E}_{\sig}(\mathcal{S}): \mathcal{N}_{r}(\mathcal{S}(\lam)) \rightarrow \mathcal{N}_{r}(\mathbb{S}_{\sigma}(\lam))$ and $\mathbb{H}_{\sig}(\mathcal{S}): \mathcal{N}_{l}(\mathcal{S}(\lam)) \rightarrow \mathcal{N}_{l}(\mathbb{S}_{\sigma}(\lam))$ are isomorphisms, then it is easy to see that $\mathbb{E}_{\sig}(G) = \mathbb{E}_{\sig}(\mathcal{S})\circ f$ and $\mathbb{H}_{\sig}(G) = \mathbb{H}_{\sig}(\mathcal{S})\circ g$. Hence $\mathbb{E}_{\sig}(\mathcal{S})\circ f$ and $\mathbb{H}_{\sig}(\mathcal{S})\circ g$ are isomorphisms. Thus $\mathbb{E}_{\sig}(G) : \mathcal{N}_{r}(G(\lam))\rightarrow \mathcal{N}_{r}(\mathbb{S}_{\sigma}(\lam))$ and $\mathbb{H}_{\sig}(\mathcal{S}) : \mathcal{N}_{l}(G(\lam)) \rightarrow \mathcal{N}_{l}(\mathbb{S}_{\sig}(\lam))$ are isomorphisms.
\end{proof}

Consider the second companion forms of $G(\lam)$ given in Example~\ref{comp}. 
Then we have 
if $\mathcal{C}_2(\lam)\left[ \begin{array}{c} v \\  \hline w \\ \end{array} \right]  = 0$ and $\left[ \begin{array}{c} u \\ \hline z \\ \end{array} \right]^{T}\mathcal{C}_2(\lam) = 0$ then  
it is easy to see that $$\left[ \begin{array}{c} v \\ \hline w \\ \end{array} \right] {\small = \left[
            \begin{array}{c}
              x \\
              (\lambda A_{m}+A_{m-1})x \\
              \vdots \\
              (\lam^{m-1}A_{m} + \lam^{m-2}A_{m-1}+\cdots + A_{1})x \\
              \hline
              (P(\lam))^{-1}Bx \\
            \end{array}
          \right]}\, \mbox{  and  }\,  \left[
                                        \begin{array}{c}
                                          u \\
                                          \hline
                                          z \\
                                        \end{array}
                                      \right]
           = \left[
     \begin{array}{c}
       \lambda^{m-1}y \\
       \lambda^{m-2}y \\
       \vdots \\
       y \\
       \hline
       (C(P(\lam))^{-1})^{T}y \\
        \end{array}
        \right]. $$

\begin{remark}
Next, consider a PGF linearization $\mathbb{K}_{\omega}(\lam)$ of the system matrix $\mathcal{S}(\lam)$ in (\ref{smh}). The recovery of zero directions of $\mathcal{S}(\lam)$ are directly follows from the Theorem~5.2 given in \cite{rafinami3}.
\end{remark}

\section{Conclusions}
we have introduced linearizations of higher order state space system. Further, we have shown that the linearized systems so obtained are {\em strict system equivalent } to the higher order systems and hence preserve system characteristics of the original systems. Furthermore, We have considered Fiedler pencils of higher order state space system and have studied recovery of zero directions of $\Sigma_1$ from those of the Fiedler-like pencils of $G(\lam).$ That is, the zero directions of the transfer functions are recovered from the eigenvectors of the Fiedler pencils without performing any arithmetic operations.

\end{document}